\documentclass{icmart} 
\usepackage{amsfonts,latexsym,color}
\input{psfig.sty}

\theoremstyle{definition}

\newtheorem{theorem}{Theorem}

\newtheorem{lemma}[theorem]{Lemma}
\newtheorem{cor}[theorem]{Corollary}
\newtheorem{conj}[theorem]{Conjecture}
\newtheorem{example}{Example}

\contact[rstan@math.mit.edu]{Richard Stanley, M.I.T.}

\begin{document}
\newcommand{\beq}{\begin{equation}} 
\newcommand{\eeq}{\end{equation}}
\newcommand{\zz}{\mathbb{Z}}
\newcommand{\pp}{\mathbb{P}} 
\newcommand{\nn}{\mathbb{N}}
\newcommand{\rr}{\mathbb{R}}
\newcommand{\cc}{\mathbb{C}}
\newcommand{\qq}{\mathbb{Q}}
\newcommand{\bm}[1]{{\mbox{\boldmath $#1$}}}
\newcommand{\con}{\mathrm{Comp}(n)}
\newcommand{\sn}{\mathfrak{S}_n}
\newcommand{\fin}{\mathfrak{I}_n}
\newcommand{\mfi}{\mathfrak{I}} 
\newcommand{\fm}{\mathfrak{M}} 
\newcommand{\fs}{\mathfrak{S}}
\newcommand{\cro}{\mathrm{cr}}
\newcommand{\nes}{\mathrm{ne}}
\newcommand{\st}{\,:\,} 
\newcommand{\ds}{\operatorname{ds}}
\newcommand{\is}{\operatorname{is}}
\newcommand{\as}{\operatorname{as}}
\newcommand{\sh}{\operatorname{sh}}
\newcommand{\rsk}{\stackrel{\mathrm{rsk}}{\rightarrow}}
\newcommand{\tworow}[2]{\genfrac{}{}{0cm}{}{#1}{#2}}

\thispagestyle{empty}

\vskip 20pt
\begin{center}
{\large\bf \textcolor{red}{Increasing and Decreasing Subsequences of
    Permutations and Their Variants}}
\vskip 15pt
{\bf \textcolor{blue}{Richard P. Stanley}}\\
{\it Department of Mathematics, Massachusetts Institute of
Technology}\\
{\it Cambridge, MA 02139, USA}\\
{\texttt{rstan@math.mit.edu}}\\[.2in]
{\bf\small version of 29 November 2005}\\
\end{center}


\abstract{We survey the theory of increasing and decreasing
  subsequences of permutations. Enumeration problems in this area are
  closely related to the RSK algorithm. The asymptotic behavior of
  the expected value of the length $\is(w)$ of the longest increasing
  subsequence of a permutation $w$ of $1,2,\dots,n$ was obtained by
  Vershik-Kerov and (almost) by Logan-Shepp. The entire limiting
  distribution of $\is(w)$ was then determined by Baik, Deift, and
  Johansson. These techniques can be applied to other classes of
  permutations, such as involutions, and are related to the
  distribution of eigenvalues of elements of the classical groups. A
  number of generalizations and variations of increasing/decreasing
  subsequences are discussed, including the theory of pattern
  avoidance, unimodal and alternating subsequences, and crossings and
  nestings of matchings and set partitions.}

\section{Introduction}
Let $\sn$ denote the symmetric group of all permutations of
$[n]:=\{1,2,\dots,n\}$. We write permutations $w\in\sn$ as
\emph{words}, i.e., $w=a_1 a_2\cdots a_n$, where $w(i)=a_i$. An
\emph{increasing subsequence} of $w$ is a subsequence $a_{i_1}\cdots
a_{i_k}$ satisfying $a_{i_1}<\cdots<a_{i_k}$, and similarly for
\emph{decreasing subsequence}. For instance, if $w=5642713$, then 567
is an increasing subsequence and 543 is a decreasing subsequence. Let
$\is(w)$ (respectively, $\ds(w)$) denote the length (number of terms)
of the longest increasing (respectively, decreasing) subsequence of
$w$. If $w=5642713$ as above, then $\is(w)=3$ (corresponding to 567)
and $\ds(w)=4$ (corresponding to 5421 or 6421). A nice interpretation
of increasing subsequences in terms of the one-person card game
\emph{patience sorting} is given by Aldous and Diaconis
\cite{a-d}. Further work on patience sorting was undertaken by
Burstein and Lankham \cite{b-l1}\cite{b-l2}. Connections between
patience sorting and airplane boarding times were found by Bachmat et
al.\ \cite{bach1}\cite{bach2} and between patience sorting and disk
scheduling by Bachmat \cite{bach3}.

The subject of increasing and decreasing subsequences began in 1935,
and there has been much recent activity. There have been major
breakthroughs in understanding the distribution of $\is(w)$, $\ds(w)$,
and related statistics on permutations, and many unexpected and deep
connections have been obtained with such areas as representation
theory and random matrix theory. A number of excellent survey papers
have been written on various aspects of these developments, e.g.,
\cite{a-d}\cite{deift}\cite{johan3}\cite{t-w3}\cite{wilf}; the
present paper will emphasize the connections with combinatorics.

In Section~\ref{sec:enum} we give some basic enumerative results
related to increasing/decreasing subsequences and show their
connection with the RSK algorithm from algebraic combinatorics. The
next two sections are devoted to the distribution of $\is(w)$ for
$w\in\sn$, a problem first raised by Ulam. In Section~\ref{sec:en} we
deal with the expectation of $\is(w)$ for $w\in\sn$, culminating in
the asymptotic formula of Logan-Shepp and Vershik-Kerov. We turn to
the entire limiting distribution of $\is(w)$ in
Section~\ref{sec:dist}. The main result is the determination of this
limiting distribution by Baik, Deift, and Johansson to be a (suitably
scaled) Tracy-Widom distribution. The Tracy-Widom distribution
originally arose in the theory of random matrices, so the result of
Baik et al.\ opens up unexpected connections between
increasing/decreasing subsequences and random matrices.

Much of the theory of increasing/decreasing subsequences of
permutations in $\sn$ carries over to permutations in certain subsets
of $\sn$, such as the set of involutions. This topic is discussed in
Section~\ref{sec:sym}. In particular, analogues of the
Baik-Deift-Johansson theorem were given by Baik and Rains. In
Section~\ref{sec3} we explain how the previous results are related to
the distribution of eigenvalues in matrices belonging to the classical
groups. 

The remaining three sections are concerned with analogues and
extensions of increasing/decreasing subsequences of
permutations. Section~\ref{sec:pav} deals with pattern avoidance,
where increasing/decreasing subsequence are replaced with other
``patterns.'' In Section~\ref{sec:alt} we consider unimodal and
alternating subsequences of permutations, and in
Section~\ref{sec:match} we replace permutations with (complete)
matchings. For matchings the role of increasing and decreasing
subsequences is played by crossings and nestings.

\textsc{Acknowledgment.} I am grateful to Percy Deift, Persi Diaconis,
Craig Tracy and Herb Wilf for providing some pertinent references.

\section{Enumeration and the RSK algorithm} \label{sec:enum}
The first result on increasing and decreasing subsequences is a famous
theorem of Erd\H{o}s and Szekeres \cite{e-s}.

\begin{theorem} \label{thm:e-s}
Let $p,q\geq 1$. If $w\in\fs_{pq+1}$, then either $\is(w)>p$ or
$\ds(w)>q$. 
\end{theorem}

This result arose in the context of the problem of determining the
least integer $f(n)$ so that any $f(n)$ points in general position in
the plane contain an $n$-element subset $S$ in convex position (i.e.,
every element of $S$ is a vertex of the convex hull of $S$). A recent
survey of this problem was given by Morris and Soltan \cite{m-s}.
Seidenberg \cite{seid} gave an exceptionally elegant proof of
Theorem~\ref{thm:e-s} based on the pigeonhole principle which has been
reproduced many times, e.g., Gardner \cite[Ch.~11, {\S}7]{gardner}.

Theorem~\ref{thm:e-s} is best possible in that there exists
$w\in\fs_{pq}$ with $\is(w)=p$ and $\ds(w)=q$. Schensted
\cite{schensted} found a quantitative strengthening of this result
based on his rediscovery of an algorithm of Robinson \cite{robinson}
which has subsequently become a central algorithm in algebraic
combinatorics. To describe Schensted's result, let
$\lambda=(\lambda_1,\lambda_2,\dots)$ be a partition of $n\geq 0$,
denoted $\lambda\vdash n$ or $|\lambda|=n$. Thus $\lambda_1\geq
\lambda_2\geq\cdots\geq 0$ and $\sum \lambda_i=n$. The \emph{(Young)
diagram} of a partition $\lambda$ is a left-justified array of
squares with $\lambda_i$ squares in the $i$th row. For instance, the
Young diagram of $(3,2,2)$ is given by
\bigskip

\centerline{\psfig{figure=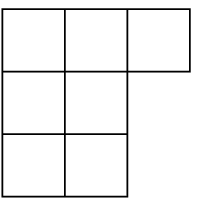}}

\bigskip
A \emph{standard Young tableau} (SYT) of shape $\lambda\vdash n$ is
obtained by placing the integers $1,2,...,n$ (each appearing once)
into the squares of the diagram of $\lambda$ (with one integer in each
square) such that every row and column is increasing. For
example, an SYT of shape $(3,2,2)$ is given by
  \[ \begin{array}{l}1\,3\,6\\ 2\,5\\4\, 7\end{array}. \]
\indent Let $f^\lambda$ denote the number of SYT of shape
$\lambda$. The quantity $f^\lambda$ has a number of additional
combinatorial interpretations \cite[Prop.~7.10.3]{ec2}. It also has a
fundamental algebraic interpretation which suggests (via
equation~(\ref{eq:gpqn}) below) a close connection between representation
theory and increasing/decreasing subsequences. Namely, the (complex)
irreducible representations $F^\lambda$ of $\sn$ are indexed in a
natural way by partitions $\lambda\vdash n$, and then 
  \beq f^\lambda = \dim F^\lambda. \label{eq:fla} \eeq
In particular (by elementary representation theory),
  \beq \sum_{\lambda\vdash n} \left( f^\lambda\right)^2 = n!. 
    \label{eq:sumfl2} \eeq
See Section~\ref{sec3} for more on the connections between
increasing/decreasing subsequences and representation theory.
MacMahon \cite[p.~175]{mac} was the first
to give a formula for $f^\lambda$ (in the guise of counting ``lattice
permutations'' rather than SYT), and later Frame, Robinson, and Thrall
\cite{f-r-t} simplified MacMahon's formula, as follows. Let $u$ be a
square of the diagram of $\lambda$, denoted $u\in\lambda$. The
\emph{hook length} $h(u)$ of (or at) $u$ is the number of squares
directly to the right or directly below $u$, counting $u$ itself once.
For instance, if $\lambda=(3,2,2)$ then the hook lengths are given by
   \[ \begin{array}{l} 5\,4\,1\\ 3\,2\\ 2\, 1 \end{array}. \]
The \emph{hook-length formula} of Frame, Robinson, and Thrall asserts
that if $\lambda\vdash n$, then 
  \beq f^\lambda = \frac{n!}{\prod_{u\in\lambda}h(u)}. 
      \label{eq:frt} \eeq
For instance, 
  \[ f^{(3,2,2)} = \frac{7!}{5\cdot 4\cdot 1\cdot 3\cdot 2\cdot 2
     \cdot 1} = 21. \]

The \emph{RSK algorithm} gives a bijection between $\sn$ and
pairs $(P,Q)$ of SYT of the same shape $\lambda\vdash n$. This
algorithm is named after Gilbert de Beauregard Robinson, who described
it in a rather vague form \cite[{\S}5]{robinson} (subsequently analyzed
by van Leeuwen \cite[{\S}7]{vl}), Craige Schensted \cite{schensted},
and Donald Knuth \cite{knuth}. For further historical information see
\cite[Notes to Ch.~7]{ec2}. The basic operation of the RSK algorithm
is \emph{row insertion}, i.e., inserting an integer $i$ into a tableau
$T$ with distinct entries and with increasing rows and columns. (Thus
$T$ satisfies the conditions of an SYT except that its entries can be
any distinct integers, not just $1,2,\dots,n$.) The process of
row inserting $i$ into $T$ produces another tableau, denoted
$T\leftarrow i$, with increasing rows and columns. If $S$ is the set
of entries of $T$, then $S\cup\{i\}$ is the set of entries of
$T\leftarrow i$. We define $T\leftarrow i$ recursively as follows. 
  \begin{itemize}
   \item If the first row of $T$ is empty or the largest entry of the
   first row of $T$ is less than 
   $i$, then insert $i$ at the end of the first row.
  \item Otherwise, $i$ replaces (or \emph{bumps}) the smallest element
  $j$ in the first row satisfying $j>i$. We then insert
  $j$ into the the second row of $T$ by the same procedure.
  \end{itemize}
For further details concerning the definition and basic properties of
$T\leftarrow i$ see e.g.\ \cite[Ch.~3]{sagan}\cite[{\S}7.11]{ec2}.

Let $w=a_1 a_2\cdots a_n\in\sn$, and let $\emptyset$ denote the empty
tableau. Define 
  \[ P_i=P_i(w) = (\cdots((\emptyset\leftarrow a_1)\leftarrow a_2)
      \leftarrow\cdots \leftarrow a_i. \]
That is, we start with the empty tableau and successively row insert
$a_1, a_2, \dots, a_i$. Set $P=P(w)=P_n(w)$. Define $Q_0=\emptyset$,
and once $Q_{i-1}$ is defined let $Q_i=Q_i(w)$ be obtained from
$Q_{i-1}$ by inserting $i$ (without changing the position of any of
the entries of $Q_{i-1}$) so that $Q_i$ and $P_i$ have the same shape.
Set $Q=Q(w)=Q_n(w)$, and finally define the output of the RSK
algorithm applied to $w$ to be the pair $(P,Q)$, denoted $w\rsk (P,Q)$.
For instance, if $w=31542\in\fs_5$, then we obtain
  \[ P_1(w) = 1,\quad P_2(w) = \begin{array}{l}1\\[-.02in] 3
  \end{array},\quad 
     P_3(w) = \begin{array}{l} 1\, 5\\[-.02in] 3 \end{array} \]
  \[ P_4(w) = \begin{array}{l} 1\, 4\\[-.02in] 3\, 5 \end{array},\quad
     P=P_5(w)=\begin{array}{l} 1\,2\\[-.02in] 3\,4\\[-.02in] 5
     \end{array}. \] 
It follows that
  \[ Q = \begin{array}{l}1\, 3\\[-.02in] 2\, 4\\[-.02in] 5 \end{array}. \]

\textsc{Note.} By a theorem of Sch\"utzenberger
\cite{schu}\cite[{\S}7.13]{ec2} we have 
  \beq Q(w)=P(w^{-1}), \label{eq:rski} \eeq 
so we could have in fact taken this formula as the definition of
$Q(w)$.

If $w\rsk (P,Q)$ and $P,Q$ have shape $\lambda$, then we also call
$\lambda$ the \emph{shape} of $w$, denoted $\lambda=\sh(w)$. The
\emph{conjugate} partition $\lambda'=(\lambda'_1,\lambda'_2,\dots)$ of
$\lambda$ is the partition whose diagram is the transpose of the
diagram of $\lambda$. Equivalently, $j$ occurs exactly
$\lambda_j-\lambda_{j+1}$ times as a part of $\lambda'$. The
\emph{length} $\ell(\lambda)$ is the number of (nonzero) parts of
$\lambda$, so $\ell(\lambda)=\lambda'_1$. The fundamental result of
Schensted \cite{schensted} connecting RSK with increasing and
decreasing subsequences is the following.

\begin{theorem} \label{thm:sch}
Let $w\in \sn$, and suppose that $\sh(w)=\lambda$. Then
  \begin{align} \is(w) & = \lambda_1 \label{eq:isw}\\
   \ds(w) & = \lambda'_1. \label{eq:dsw} \end{align}
\end{theorem}

Equation (\ref{eq:isw}) is easy to prove by induction since we
need only analyze the effect of the RSK algorithm on the first row of
the $P_i$'s. On the other hand, equation
(\ref{eq:dsw}) is based on the following symmetry
property of RSK proved by Schensted. If $w=a_1 a_2\cdots a_n$ then let
$w^r=a_n \cdots a_2 a_1$, the \emph{reverse} of $w$. We then have
  \beq w\rsk (P,Q)\ \Rightarrow\ w^r\rsk (P^t,\mathrm{evac}(Q)^t), 
   \label{eq:wr} \eeq
where $^t$ denotes transpose and evac$(Q)$ is a certain tableau
called the \emph{evacuation} of $Q$ (first defined by Sch\"utzenberger
\cite{schu2}) which we will not define here. Equation~(\ref{eq:wr})
shows that if $\sh(w)=\lambda$, then $\sh(w^r)=\lambda'$. Since
clearly $\is(w^r)=\ds(w)$, equation (\ref{eq:dsw}) follows from
(\ref{eq:isw}).  

Theorem~\ref{thm:sch} has several immediate consequences. The first
corollary is the Erd\H{o}s-Szekeres theorem (Theorem~\ref{thm:e-s}),
for if $\sh(w)=\lambda$, $\is(w)\leq p$, and $\ds(w)\leq q$, then
$\lambda_1\leq p$ and $\lambda'_1\leq q$. Thus the diagram of
$\lambda$ is contained in a $q\times p$ rectangle, so $|\lambda|\leq
pq$. By the same token we get a quantitative statement that
Theorem~\ref{thm:e-s} is best possible.

\begin{cor} \label{cor:pq}
The number of permutations $w\in \fs_{pq}$, where say $p\leq q$,
satisfying $\is(w)=p$ and $\ds(w)=q$ is given by
  \beq (f^{(p^q)})^2 = \left(\frac{(pq)!}{1^12^2\cdots p^p(p+1)^p
    \cdots q^p(q+1)^{p-1}\cdots (p+q-1)^1}\right)^2, 
    \label{eq:pq} \eeq
where $(p^q)$ denotes the partition with $q$ parts equal to $p$.
\end{cor} 

\begin{proof}
Let $\lambda=\sh(w)$. If $\is(w)=p$ and $\ds(w)=q$, then $\lambda_1=p$
and $\lambda'_1=q$. Since $\lambda\vdash pq$, we must have
$\lambda=(p^q)$. The number of $v\in\sn$ with a fixed shape $\mu$ is
just $(f^\mu)^2$, the number of pairs $(P,Q)$ of SYT of shape
$\mu$. Hence the left-hand side of equation (\ref{eq:pq}) follows. The
right-hand side is then a consequence of the hook-length formula
(\ref{eq:frt}). 
\end{proof}

\begin{figure}
\centering
 \centerline{\psfig{figure=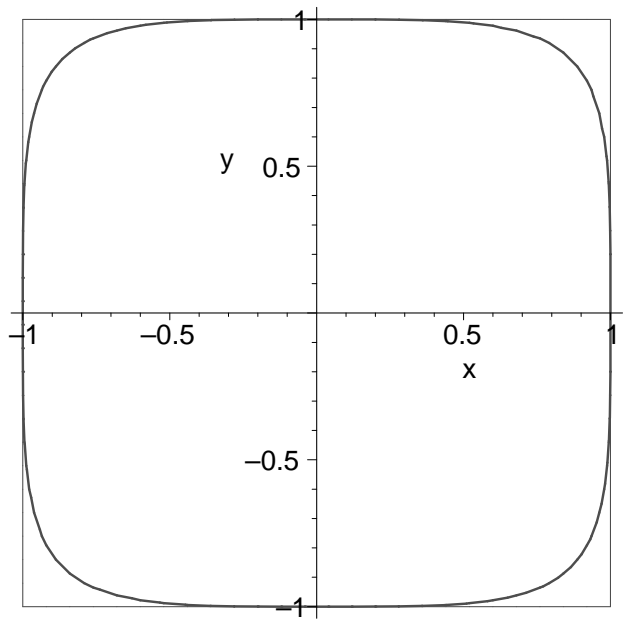}}
\caption{The curve $(x^2-y^2)^2+2(x^2+y^2)= 3$}
\label{fig:romik}
\end{figure}

An interesting result concerning the extremal permutations in the case
$p=q$ in Corollary~\ref{cor:pq} was obtained by Romik
\cite[Thm.~5]{romik}. It can be stated informally as follows. Pick a
random permutation $w\in\fs_{p^2}$ satisfying $\is(w)=\ds(w)=p$. Let
$P_w$ be the $p^2\times p^2$ permutation matrix corresponding to $w$,
drawn in the plane so that its corners occupy the points $(\pm 1,\pm
1)$. Then almost surely as $p\rightarrow\infty$ the limiting curve
enclosing most of the 1's in $P_w$ is given by
  \[ \{(x,y)\in\rr^2\st (x^2-y^2)^2+2(x^2+y^2)= 3\}. \] 
See Figure~\ref{fig:romik}. In particular, this curve encloses a
fraction $\alpha=0.94545962\cdots$ of the entire square with vertices
$(\pm 1,\pm 1)$. The number $\alpha$ can be expressed in terms of
elliptic integrals as
  \[ \alpha = 2\int_0^1\frac{1}{\sqrt{(1-t^2)(1-(t/3)^2)}}dt
    -\frac 32\int_0^1\sqrt{\frac{1-(t/3)^2}{1-t^2}}dt. \]
Compare with the case of \emph{any} $w\in \fs_n$, when
clearly the limiting curve encloses the entire square with vertices
$(\pm 1,\pm 1)$. For further information related to permutations
  $w\in\fs_{p^2}$ satisfying $\is(w)=\ds(w)=p$, see the paper
  \cite{p-r} of Pittel and Romik.

Clearly Corollary~\ref{cor:pq} can be extended to give a formula
\cite[Cor.~7.23.18]{ec2} for the number 
  \[ g_{pq}(n) = \#\{w\in\sn\st \is(w)=p,\ \ds(w)=q\}, \]
namely,
  \beq g_{pq}(n) = \sum_{\substack{\lambda\vdash n\\ \lambda_1=p,\
  \lambda'_1= q}} \left( f^\lambda\right)^2. \label{eq:gpqn} \eeq
The usefulness of this formula may not be readily apparent, but
Theorem~\ref{thm:regev} below is an example of its utility.

\section{Expectation of is$(w)$} \label{sec:en}
A further application of Theorem~\ref{thm:sch} concerns the
distribution of the function $\is(w)$ as $w$ ranges over $\sn$. The
problem of obtaining information on this distribution was first raised
by Ulam \cite[{\S}11.3]{ulam} in the context of Monte Carlo calculations. 
In particular, one can ask for
information on the expectation $E(n)$ of $\is(w)$ for $w\in\sn$, i.e., 
  \[ E(n) = \frac{1}{n!}\sum_{w\in\sn} \is(w). \]
Ulam mentions the computations of E. Neighbor suggesting that $E(n)$ is
about $1.7\sqrt{n}$. Numerical experiments by Baer and Brock
\cite{b-b} suggested that $E(n)\sim 2\sqrt{n}$ might be closer to the
truth. The Erd\H{o}s-Szekeres theorem (Theorem~\ref{thm:e-s}) implies
immediately that $E(n)\geq \sqrt{n}$, since
  \[ \frac 12(\is(w)+\is(w^r))\geq \sqrt{\is(w)\is(w^r)} =
      \sqrt{\is(w)\ds(w)}\geq \sqrt{n}. \]
Hammersley \cite{hamm} was the first person to seriously consider the
question of estimating $E(n)$. He showed that if 
  \[ c=\lim_{n\rightarrow\infty} \frac{E(n)}{\sqrt{n}}, \]
then $c$ exists and satisfies 
    \[ \frac{\pi}{2} \leq c \leq e. \]
He also gave a heuristic argument that $c=2$, in agreement with the
experiments of Baer and Brock.

The next progress on Ulam's problem was based on Schensted's theorem
(Theorem~\ref{thm:sch}). It follows from this result that
  \beq E(n) = \frac{1}{n!}\sum_{\lambda\vdash n} \lambda_1\left(
     f^\lambda\right)^2. \label{eq:en} \eeq
Now the RSK algorithm itself shows that
  \beq n! = \sum_{\lambda\vdash n} \left( f^\lambda\right)^2,
    \label{eq:fl2} \eeq
in agreement with (\ref{eq:sumfl2}).
Since the number of terms in the sum on the right-hand side of
(\ref{eq:fl2}) is very small compared to $n!$, the maximum value of
$f^\lambda$ for $\lambda\vdash n$ is close to $\sqrt{n!}$. Let
$\lambda^n$ be a value of $\lambda\vdash n$ for which $f^\lambda$ is
maximized. Then by (\ref{eq:en}) we see that a close approximation to
$E(n)$ is given by
 \begin{align*} E(n) & \approx \frac{1}{n!}(\lambda^n)_1\left(
    f^{\lambda^n}\right)^2\\  & \approx (\lambda^n)_1. \end{align*}
This heuristic argument shows the importance of determining the
partition $\lambda^n$ maximizing the value of $f^\lambda$ for
$\lambda\vdash n$. 

We are only really interested in the behavior of $\lambda^n$ as
$n\rightarrow\infty$, so let us normalize the Young diagram of any
partition $\lambda$ to have area one. Thus each square of the diagram
has length $1/\sqrt{n}$. Let the upper boundary of (the
diagram of) $\lambda$ be the $y$-axis directed to the right, and the
left boundary be the $x$-axis directed down. As $n\rightarrow \infty$
it is not unreasonable to expect that the boundary of the partition
$\lambda^n$ will approach some limiting curve $y=\Psi(x)$. If this
curve intersects the $x$-axis at $x=b$, then it is immediate that
  \[ c:=\lim_{n\rightarrow\infty} \frac{E(n)}{\sqrt{n}}\geq b. \]
We cannot be sure that $b=c$ since conceivably the first few parts of
$\lambda^n$ are much larger than the other parts, so these parts would
``stretch out'' the curve $y=\Psi(x)$ along the $x$-axis.

\begin{figure}
\centering
 \centerline{\psfig{figure=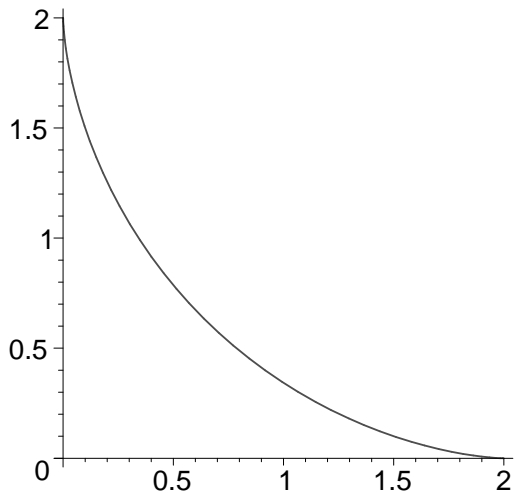}}
\caption{The curve $y=\Psi(x)$}
\label{fig:psi}
\end{figure}

It was shown independently by Vershik-Kerov \cite{v-k} and Logan-Shepp
\cite{l-s}  that $y=\Psi(x)$ indeed does exist and is given
parametrically by 
  \begin{align*} x & =  y + 2\cos \theta\\
    y & =  \frac{2}{\pi}(\sin\theta-\theta\cos\theta), \end{align*}
for $0\leq \theta\leq \pi$. See Figure~\ref{fig:psi}, where we have
placed the coordinate axes in their customary locations. 

Logan-Shepp and Vershik-Kerov obtain the curve $\Psi(x)$ as a solution
to a variational problem. If $(x,y)$ is a point in the region $A$
enclosed by the curve and the coordinate axes, then the (normalized)
hook-length at $(x,y)$ is $f(x)-y+f^{-1}(y)-x$. By equation
(\ref{eq:frt}) we maximize $f^\lambda$ by minimizing
$\prod_{u\in\lambda}h(u)$. Hence in the limit we wish to minimize
$I(f)=\int\!\!\int_A \log(f(x)-y+f^{-1}(y)-x)dx\,dy$, subject to the
normalization $\int\!\!\int_A dx\,dy=1$. It is shown in \cite{l-s} and
\cite{v-k} that $f=\Psi$ is the unique function minimizing $I(f)$ (and
moreover $I(f)=-1/2$).

\textsc{Note.} If we extend the curve $y=\Psi(x)$ to include the
$x$-axis for $x\geq 2$ and the $y$-axis for $y\geq 2$, and then rotate
it 45$^\circ$ counterclockwise, then we obtain the curve
  \[ \Omega(x) = \left\{ \begin{array}{rl}
   \frac{2}{\pi}(x\arcsin(x/2)+\sqrt{4-x^2}), & |x|\leq 2,\\[.05in]
   |x|, & |x|\geq 2. \end{array} \right. \]
This form of the limiting curve is more convenient for some purposes
\cite{kerov}\cite{kerov2}\cite{kerov3}, such as surprising connections
with the separation of zeros of orthogonal polynomials. 

We see immediately from the equations for $\Psi(x)$ that it intersects
the $x$-axis at $x=2$, so $c\geq 2$. By a simple but clever use of the
RSK algorithm Vershik and Kerov show in their paper that $c\leq 2$, so
we conclude that 
  \beq E(n)\sim 2\sqrt{n}. \label{eq:vkexp} \eeq 


Different proofs that $E(n)\sim 2\sqrt{n}$ were later given by Aldous
and Diaconis \cite{a-d}, Groeneboom \cite{groen}, Johansson
\cite{johan2}, and Sep\"al\"ainen 
\cite{sep}.  The proof of Aldous and Diaconis is based on an
interacting particle process for which the number of particles
remaining after $n$ steps has the same distribution as is$_n$. Their
proof is known in the language of statistical physics as a
\emph{hydrodynamic limit} argument. See \cite[{\S}3]{a-d2} for a brief
survey.

We should remark that the curve $y=\Psi(x)$ is not merely the limiting
curve for the partition \emph{maximizing} $f^\lambda$; it is also the
limiting curve for the \emph{typical} shape of a permutation
$w\in\sn$. A remarkable refinement of this fact is due to Kerov
\cite{i-o}\cite{kerov:gp}\cite[{\S}0.3.4]{kerov3}, who shows that the
deviation of a Young diagram from the expected limit converges in
probability to a certain Gaussian process. A different kind of
refinement is due to Borodin, Okounkov, and Olshanski
\cite[Thm.~1]{b-o-o}, who obtain more detailed local information about
a typical shape $\lambda$ than is given by $\Psi(x)$.

\section{Distribution of is$(w)$} \label{sec:dist}
A major breakthrough in understanding the behavior of $\is(w)$ was
achieved in 1999 by Baik, Deift, and Johansson \cite{b-d-j}. They
determined the entire limiting distribution of $\is(w)$ as
$n\rightarrow\infty$. It turns out to be given by the (suitably
scaled) Tracy-Widom distribution, which had first appeared in
connection with the distribution of the largest eigenvalue of a random
hermitian matrix.

To describe these results, write $\is_n$ for the function
$\is:\sn\rightarrow\zz$. 
Let $u(x)$ denote the unique solution to the nonlinear second order
differential equation
  \beq u''(x) = 2u(x)^3+xu(x), \label{eq:pain} \eeq
subject to the condition
 \[ u(x) \sim -\frac{e^{-\frac 23 x^{3/2}}}
    {2\sqrt{\pi}x^{1/4}},\ \mathrm{as}\ x\rightarrow\infty. \]
Equation (\ref{eq:pain}) is known as the \emph{Painlev\'e II
  equation}, after Paul Painlev\'e (1863--1933). Painlev\'e
completely classified differential equations (from a certain class of
second order equations) whose ``bad'' singularities (branch points and
essential singularities) were independent of the initial conditions.
Most of the equations in this class were already known, but a few were
new, including equation (\ref{eq:pain}).

Now define the \emph{Tracy-Widom distribution} to be the probability
distribution on $\rr$ given by 
  \beq F(t)= \exp\left( -\int_t^\infty (x-t)u(x)^2\,dx\right).
   \label{eq:t-w} \eeq
It is easily seen that $F(t)$ is indeed a probability distribution,
i.e., $F(t)\geq 0$ and $\int_{-\infty}^\infty F'(t)dt=1$. 
We can now state the remarkable results of Baik, Deift, and
Johansson. 

\begin{theorem} \label{thm:bdj}
We have for random (uniform) $w\in\sn$ and all $t\in\rr$ that
  \[ \lim_{n\rightarrow\infty} \mathrm{Prob}\left(
    \frac{\is_n(w)-2\sqrt{n}}{n^{1/6}}\leq t\right) =F(t). \]  
\end{theorem}

The above theorem is a vast refinement of the Vershik-Kerov and
Logan-Shepp results concerning $E(n)$, the expectation of is$(w)$. It
gives the entire limiting distribution (as $n\rightarrow\infty$) of
is$_n(w)$. Baik, Deift, and Johansson also determine all the limiting
moments of $\is_n(w)$. In particular, we have the following formula
for the variance Var$(\is_n)$ of $\is_n$ as $n\rightarrow\infty$.

\begin{cor} \label{cor:varis}
We have
  \begin{align*} \lim_{n\rightarrow\infty}\frac{\mathrm{Var}(\is_n)}
     {n^{1/3}} & =  \int t^2\,dF(t) -\left(\int
       t\,dF(t)\right)^2\\ & = 0.8131947928\cdots, \end{align*}
and
   \begin{align} \lim_{n\rightarrow\infty}\frac{E(n)-2\sqrt{n}}
      {n^{1/6}} & =  \int t\,dF(t) \label{eq:e2} \\
     & =  -1.7710868074\cdots. \nonumber \end{align}
\end{cor}

Note that equation (\ref{eq:e2}) may be rewritten
  \[ E(n) = 2\sqrt{n}+\alpha n^{1/6}+o(n^{1/6}), \]
where $\alpha=\int t\,dF(t)$, thereby giving the second term in the
asymptotic behavior of $E(n)$. 

We will say only a brief word on the proof of Theorem~\ref{thm:bdj},
explaining how combinatorics enters into the picture. Some kind of
analytic expression is needed for the distribution of is$_n(w)$. Such
an expression is provided by the following result of Ira Gessel
\cite{gessel}, later proved in other ways by various persons; see
\cite[{\S}1]{b-d-j} for references. Define
  \begin{align} u_k(n) & = \#\{ w\in\sn\st\is_n(w)\leq k\}
   \label{eq:ukn} \\
    U_k(x) & = \sum_{n\geq 0} u_k(n)\frac{x^{2n}}{n!^2},\ k\geq 1 
     \nonumber \\ 
    I_i(2x) & = \sum_{n\geq 0}\frac{x^{2n+i}}{n!\,(n+i)!},\
       i\in\zz. \nonumber \end{align}
The function $I_i$ is the \emph{hyperbolic Bessel function} of the
first kind of order $i$. Note that $I_i(2x)=I_{-i}(2x)$.

\begin{theorem} \label{thm:gessel}
We have
  \[ U_k(x) = \det\left( I_{i-j}(2x)\right)_{i,j=1}^k. \]
\end{theorem}

\begin{example} \label{ex:ukn} 
We have (using $I_1=I_{-1}$)
  \begin{align*} U_2(x) & =  \left| \begin{array}{cc} I_0(2x) & I_1(2x)\\
                I_1(2x) & I_0(2x) \end{array} \right|\\
           & =  I_0(2x)^2-I_1(2x)^2. \end{align*}
From this expression it is easy to deduce that 
  \beq u_2(n) = \frac{1}{n+1}{2n\choose n}, \label{eq:u2n} \eeq
a Catalan number. This formula for $u_2(n)$ was first stated by 
Hammersley \cite{hamm} in 1972, with the first published proofs by
Knuth \cite[{\S}5.1.4]{knuth2} and Rotem \cite{rotem}. There is a
more complicated expression for $u_3(n)$ due to Gessel
\cite[{\S}7]{gessel}\cite[Exer.~7.16(e)]{ec2}, namely,
  \beq u_3(n) =\frac{1}{(n+1)^2(n+2)} \sum_{j=0}^n {2j\choose j}
    {n+1\choose j+1}{n+2\choose j+2}, \label{eq:u3n} \eeq
while no ``nice'' formula for $u_k(n)$ is known for fixed $k>3$. It is
known, however, that $u_k(n)$ is a P-recursive function of $n$, i.e.,
satisfies a linear recurrence with polynomial coefficients
\cite[{\S}7]{gessel}. For instance,
  \begin{align*} (n+4)(n+3)^2u_4(n) & = (20n^3+62n^2+22n-24)u_4(n-1)\\
     & \quad -64n(n-1)^2u_4(n-2)\\[.2in]
  (n+6)^2(n+4)^2u_5(n) & = (375-400n-843n^2-322n^3-35n^4)u_5(n-1)\\
    & \quad +(259n^2+622n+45)(n-1)^2u_5(n-2)\\ &
    \quad -225(n-1)^2(n-2)^2u_5(n-3). \end{align*}
A number of conjectures about the form of the recurrence satisfied by
$u_k(n)$ were made by Bergeron, Favreau, and Krob \cite{b-f-k},
reformulated in \cite{b-g} with some progress toward a proof.
\end{example}

Gessel's theorem (Theorem~\ref{thm:gessel}) reduces the theorem of
Baik, Deift, and Johansson to ``just'' analysis, viz., the
Riemann-Hilbert problem in the theory of integrable systems, followed
by the method of steepest descent to analyze the asymptotic behavior
of integrable systems. For further information see the survey
\cite{deift} of Deift mentioned above. 

The asymptotic behavior of is$_n(w)$ (suitably scaled) turned out to
be identical to the Tracy-Widom distribution $F(t)$ of equation
(\ref{eq:t-w}). Originally the Tracy-Widom distribution arose in
connection with the \emph{Gaussian Unitary Ensemble} (GUE). GUE is a
certain natural probability density on the space of all $n\times
n$ hermitian matrices $M=(M_{ij})$, namely,
  \[ Z_n^{-1} e^{-\mathrm{tr}(M^2)}dM, \]
where $Z_n$ is a normalization constant and
  \[ dM = \prod_i dM_{ii}\cdot\prod_{i<j} d(\mathrm{Re}\, M_{ij})
   d(\mathrm{Im}\,M_{ij}). \]
Let the eigenvalues of $M$ be $\alpha_1\geq \alpha_2\geq \cdots \geq
\alpha_n$. The following result marked the eponymous appearance
\cite{t-w} of the Tracy-Widom distribution:
  \beq \lim_{n\rightarrow\infty}\mathrm{Prob}\left( \left(\alpha_1-
    \sqrt{2n}\right)\sqrt{2}n^{1/6}\leq t\right) = F(t). \label{eq:tw}
   \eeq
Thus as $n\rightarrow\infty$, is$_n(w)$ and $\alpha_1$ have the same
distribution (after scaling). 

It is natural to ask, firstly, whether there is a result analogous to
equation (\ref{eq:tw}) for the other eigenvalues $\alpha_k$ of the GUE
matrix $M$, and, secondly, whether there is some connection between
such a result and the behavior of increasing subsequences of random
permutations.  A generalization of (\ref{eq:tw}) to all $\alpha_k$ was
given by Tracy and Widom \cite{t-w} (expressed in terms of the
Painlev\'e II function $u(x)$). The connection with increasing
subsequences was conjectured in \cite{b-d-j} and proved independently
by Borodin-Okounkov-Olshanski \cite{b-o-o}, Johannson \cite{johan},
and Okounkov \cite{ok}, after first being proved for the second
largest eigenvalue by Baik, Deift, and Johansson \cite{b-d-j2}. Given
$w\in\sn$, define integers 
$\lambda_1,\lambda_2,\dots$ by letting $\lambda_1+\cdots+\lambda_k$ be
the largest number of elements in the union of $k$ increasing
subsequences of $w$. For instance, let $w=247951368$. The longest
increasing subsequence is 24568, so $\lambda_1=5$. The largest union
of two increasing subsequences is 24791368 (the union of 2479 and
1368), so $\lambda_1+\lambda_2=8$.  (Note that it is impossible to
find a union of length 8 of two increasing subsequences that contains
an increasing subsequence of length $\lambda_1=5$.) Finally $w$ itself
is the union of the three increasing subsequences 2479, 1368, and 5,
so $\lambda_1+\lambda_2+\lambda_3=9$. Hence
$(\lambda_1,\lambda_2,\lambda_3)=(5,3,1)$ (and $\lambda_i=0$ for
$i>3$). Readers familiar with the theory of the RSK algorithm will
recognize the sequence $(\lambda_1, \lambda_2, \dots)$ as the shape
sh$(w)$ as defined preceding Theorem~\ref{thm:sch}, a well-known
result of Curtis Greene \cite{greene}\cite[Thm.\ A1.1.1]{ec2}. (In
particular, $\lambda_1\geq \lambda_2\geq\cdots$, a fact which is by no
means obvious.)  The result of \cite{b-o-o}\cite{johan}\cite{ok}
asserts that as as $n\rightarrow\infty$, $\lambda_k$ and $\alpha_k$
are equidistributed, up to scaling. In particular, the paper \cite{ok}
of Okounkov provides a direct connection, via the topology of random
surfaces, between the two seemingly unrelated appearances of the
Tracy-Widom distribution in the theories of random matrices and
increasing subsequences. A very brief explanation of this connection
is the following: a surface can be described either by gluing together
polygons along their edges or by a ramified covering of a sphere. The
former description is related to random matrices via the theory of
quantum gravity, while the latter can be formulated in terms of the
combinatorics of permutations.

We have discussed how Gessel's generating function $U_k(x)$ for
$u_k(n)$ is needed to find the limiting distribution of $\is_n$. We
can also ask about the behavior of $u_k(n)$ itself for fixed $k$. The
main result here is due to Regev \cite{regev}.

\begin{theorem} \label{thm:regev}
For fixed $k$ and for $n\rightarrow\infty$ we have the asymptotic
formula 
  \[ u_k(n) \sim 1!\,2!\cdots(k-1)!
    \left(\frac{1}{\sqrt{2\pi}}\right)^{k-1}
     \left(\frac 12\right)^{(k^2-1)/2}k^{k^2/2}
      \frac{k^{2n}}{n^{(k^2-1)/2}}. \]
\end{theorem}

\textbf{Idea of proof.} From the RSK algorithm we have
  \beq u_k(n) = \sum_{\tworow{\lambda\vdash n}
     {\ell(\lambda)\leq k}}  \left( f^\lambda\right)^2. 
     \label{eq:uknsum} \eeq
Write $f^\lambda$ in terms of the hook-length formula, factor out the
dominant term from the sum (which can be determined via Stirling's
formula), and interpret what remains in the limit $n\rightarrow
\infty$ as a $k$-dimensional integral. This integral turns out to be a 
special case of \emph{Selberg's integral} (e.g., \cite[Ch.~8]{a-a-r}),
which can be explicitly evaluated. $\ \Box$

An immediate corollary of Theorem~\ref{thm:regev} (which can also be
easily proved directly using RSK) is the formula
  \beq \lim_{n\rightarrow \infty} u_k(n)^{1/n} = k^2. 
     \label{eq:limuk} \eeq

\section{Symmetry.} \label{sec:sym}
Previous sections dealt with properties of general
permutations in $\sn$. Much of the theory carries over for certain
classes of permutations. There is a natural action of the dihedral
group $D_4$ of order 8 on $\sn$, best understood by considering the
permutation matrix $P_w$ corresponding to $w\in\sn$. Since $P_w$ is a
square matrix, $D_4$ acts on $P_w$ as the usual symmetry group of the
square. In particular, reflecting through the main diagonal transforms
$P_w$ to its transpose $P_w^t=P_{w^{-1}}$. Reflecting about a
horizontal line produces $P_{w^r}$, where $w^r$ is the reverse of $w$
as used in equation~(\ref{eq:wr}). These two reflections generate the
entire group $D_4$.

Let $G$ be a subgroup of $D_4$, and let 
  \[ \sn^G = \{w\in\sn\st \sigma\cdot w=w\ \mathrm{for\ all}\
     \sigma\in G\}. \] 
Most of the results of the preceding sections can be carried over from
$\sn$ to $\sn^G$.  The general theory is due to Baik and Rains
\cite{b-r3}\cite{b-r1}\cite{b-r2}. Moreover, for certain $G$ we can
add the condition that no entry of $P_w$ equal to 1 can be fixed by
$G$, or more strongly we can specify the number of 1's in $P_w$ fixed
by $G$.  For instance, if $G$ is the group of order 2 generated by
reflection through the main diagonal, then we are specifying the
number of fixed points of $w$. For convenience we will consider here
only two special cases, viz., (a) $G$ is the group of order 2
generated by reflection through the main diagonal. In this case $\sn^G
=\{w\in\sn \st w^2=1\}$, the set of \emph{involutions} in $\sn$, which
we also denote as $\fin$.  (b) The modification of (a) where we
consider fixed-point free involutions only. Write $\fin^*$ for this
set, so $\fin^*=\emptyset$ when $n$ is odd.

The RSK algorithm is well-behaved with respect to inversion, viz., it
follows from equation (\ref{eq:rski}) that if $w\rsk (P,Q)$ then
$w^{-1}\rsk (Q,P)$. Hence 
  \beq w^2=1\ \mathrm{if\ and\ only\ if}\ P=Q. \label{eq:irsk} \eeq
Let
  \[ y_k(n) = \#\{w\in \fin\st \is_n(w)\leq k\}. \] 
By Schensted's theorem (Theorem~\ref{thm:sch}) we conclude
  \[ y_k(n) = \sum_{\tworow{\lambda\vdash
    n}{\lambda_1\leq k}} f^\lambda, \]
the ``involution analogue'' of (\ref{eq:uknsum}). From this formula or
by other means one can obtain formulas for $y_k(n)$ for small $k$
analogous to (\ref{eq:u2n}) and (\ref{eq:u3n}). In particular (see
\cite[Exer.~7.16(b)]{ec2} for references),
   \begin{align*} y_2(n) & = {n\choose \lfloor n/2 \rfloor}\\
   y_3(n) & = \sum_{i=0}^{\lfloor n/2 \rfloor}{n\choose 2i}C_i\\
   y_4(n) & = C_{\lfloor (n+1)/2 \rfloor}C_{\lceil (n+1)/2 \rceil}\\
   y_5(n) & = 6\sum_{i=0}^{\lfloor n/2 \rfloor}{n\choose 2i}C_i
        \frac{(2i+2)!}{(i+2)!(i+3)!}, \end{align*}
where as usual $C_i$ is a Catalan number.

The RSK algorithm is also well-behaved with respect to fixed points of
involutions. It was first shown by Sch\"utzenberger
\cite[p.~127]{schu}\cite[Exer.~7.28(a)]{ec2} that if $w^2=1$ and
$w\rsk (P,P)$, then the number of fixed points of $w$ is
equal to the number of columns of $P$ of odd length. Let 
  \begin{align*} v_{2k}(n) & = \#\{w\in\fin^*\st \ds(w)\leq 2k\}\\
    z_k(n) & = \#\{w\in\fin^*\st \is(w)\leq k\}. \end{align*}
(It is easy to see directly that if $w\in\fin^*$ then $\ds(w)$ is
even, so there is no need to deal with $v_{2k+1}(n)$.) It follows that 
  \begin{align*} v_{2k}(n) & =  \sum_{\substack{\lambda\vdash n\\
        \lambda_1\leq k}} f^{2\lambda'}\\
     z_k(n) & =  \sum_{\substack{\lambda\vdash n\\ \lambda'_1\leq k}}
   f^{2\lambda'}, \end{align*}
where $2\lambda'=(2\lambda'_1,2\lambda'_2,\dots)$, the general
partition with no columns of odd length. Note that for fixed-point
free involutions $w\in\fin^*$ we no longer have a symmetry between
$\is(w)$ and $\ds(w)$, as we do for arbitrary permutations or
arbitrary involutions. 

There are also ``involution analogues'' of Gessel's determinant
(Theorem~\ref{thm:gessel}). Equations (\ref{eq:y2ke}) and
(\ref{eq:y2ko}) below were first obtained by Gessel
\cite[{\S}6]{gessel}, equation (\ref{eq:v2k}) by Goulden
\cite{goulden}, and equations (\ref{eq:z2ke}) and (\ref{eq:z2ko}) by
Baik and Rains \cite[Cor.~5.5]{b-r3}. Let
  \begin{align*} Y_k(x) & =  \sum_{n\geq 0}y_k(n)\frac{x^n}{n!}\\
        V_{2k}(x) & =  \sum_{n\geq 0}v_{2k}(n)\frac{x^n}{n!}\\
        Z_k(x) & =  \sum_{n\geq 0}z_k(n)\frac{x^n}{n!}. \end{align*}
Write $I_i=I_i(2x)$.

\begin{theorem} \label{thm:inv}
We have
  \begin{align} Y_{2k}(x) & =  \det\left( I_{i-j}+I_{i+j-1}
            \right)_{i,j=1}^k \label{eq:y2ke}\\ 
        Y_{2k+1}(x) & =  e^x\det\left( I_{i-j}-I_{i+j}
            \right)_{i,j=1}^k \label{eq:y2ko}\\
        V_{2k}(x) & =  \det\left( I_{i-j}-I_{i+j}
           \right)_{i,j=1}^k \label{eq:v2k}\\
        Z_{2k}(x) & =  \frac 14 \det\left( I_{i-j}+I_{i+j-2}
           \right)_{i,j=1}^k +\frac 12 
          \det\left( I_{i-j}-I_{i+j} \right)_{i,j=1}^{k-1}
                 \label{eq:z2ke}\\ 
       Z_{2k+1}(x) & =  \frac 12 e^x \det\left( I_{i-j}-I_{i+j-1}
           \right)_{i,j=1}^k +\frac 12 e^{-x}
        \det\left( I_{i-j}+I_{i+j-1}
           \right)_{i,j=1}^k. \label{eq:z2ko}
          \end{align}
\end{theorem}

Once we have the formulas of Theorem~\ref{thm:inv} we can use the
techniques of Baik, Deift, and Johansson to obtain the limiting
behavior of $\ds(w)$ for $w\in\fin$ and $w\in\fin^*$. These results
were first obtained by Baik and Rains \cite{b-r1}\cite{b-r2}.

\begin{theorem} \label{thm:sym}
\rm{(a)} We have for random (uniform) $w\in\fin$ and all $t\in\rr$
that 
  \[ \lim_{n\rightarrow\infty} \mathrm{Prob}
  \left(\frac{\is_n(w)-2\sqrt{n}}{n^{1/6}}\leq t\right) 
     = F(t)^{1/2}\exp\left( \frac 12\int_t^\infty u(s)ds\right), \]
where $F(t)$ denotes the Tracy-Widom distribution and $u(s)$ the
Painlev\'e II function. (By \eqref{eq:irsk} we can replace
$\is_n(w)$ with $\ds_n(w)$.)\\ \indent
\rm{(b)} We have for random (uniform) $w\in\mfi_{2n}^*$ and all
  $t\in\rr$ that 
  \[ \lim_{n\rightarrow\infty} \mathrm{Prob}
  \left(\frac{\ds_{2n}(w)-2\sqrt{2n}}{(2n)^{1/6}}\leq t\right) 
     = F(t)^{1/2}\exp\left( \frac 12\int_t^\infty u(s)ds\right). \]
\indent
  \rm{(c)} We have for random (uniform) $w\in\mfi_{2n}^*$ and all
  $t\in\rr$ that 
  \[ \lim_{n\rightarrow\infty} \mathrm{Prob}
  \left(\frac{\is_{2n}(w)-2\sqrt{2n}}{(2n)^{1/6}}\leq t\right) 
     =  F(t)^{1/2}\cosh\left( \frac 12\int_t^\infty u(s)ds\right). \]
\end{theorem}

There are orthogonal and symplectic analogues of the GUE model of
random hermitian matrices, known as the GOE and GSE models. The GOE
model replaces hermitian matrices with real symmetric matrices, while
the GSE model concerns hermitian self-dual matrices. (A $2n\times 2n$
complex matrix is \emph{hermitian self-dual} if is composed of
$2\times 2$ blocks of the form $\left[ \begin{array}{ll} a+ bi &
c+di\\ -c+di & a-bi \end{array}\right]$ which we identify with the
quaternion $a+bi+cj+dk$, such that if we regard the matrix as an
$n\times n$ matrix $M$ of quaternions, then $\overline{M}_{ji}=M_{ij}$
where the bar is quaternion conjugation.)
The limiting distribution of $\ds_{2n}(w)$ for $w\in\mfi_{2n}^*$
coincides (after scaling) with the distribution of the largest
eigenvalue of a random real symmetric matrix (under the GOE model),
while the limiting distribution of $\is_{2n}(w)$ for $w\in\mfi_{2n}^*$
coincides (after scaling) with the distribution of the largest
eigenvalue of a random hermitian self-dual matrix (under the GSE
model) \cite{t-w2}.

\section{Connections with the classical groups.} \label{sec3}
In equation (\ref{eq:gpqn}) we expressed $g_{pq}(n)$, the number of
permutations $w\in\sn$ satisfying $\is(w)=p$ and $\ds(w)=q$, in terms
of the degrees $f^\lambda$ of irreducible representations of $\sn$. 
This result can be restated via Schur-Weyl duality as a statement
about the distribution of eigenvalues of matrices in the unitary group
$U(n)$. The results of Section~\ref{sec:sym} can be used to extend
this statement to other classical groups. 

Let $U(k)$ denote the group of $k\times k$ complex unitary
matrices. For a function $f:U(k)\rightarrow \cc$, let $E(f)$
denote expectation with respect to Haar measure, i.e., 
  \[ E(f) = \int_{M\in U(k)} f(M)dM, \]
where $\int$ is the Haar integral. The following result was proved by
Diaconis and Shahshahani \cite{d-s} for $n\geq k$ and by Rains
\cite{rains1} for general $k$. Note that if $M$ has
eigenvalues $\theta_1,\dots,\theta_k$ then
  \[ |\mathrm{tr}(M)^n|^2 = (\theta_1+\cdots+\theta_k)^n
          (\bar{\theta}_1+\cdots+\bar{\theta}_k)^n. \]

\begin{theorem} \label{thm:trmn}
We have $E(|\mathrm{tr}(M)^n|^2) = u_k(n)$, where $u_k(n)$ is defined
in equation (\ref{eq:ukn}).
\end{theorem}

\begin{proof}
The proof is based on the theory of symmetric functions, as developed
e.g.\ in \cite{macd} or \cite[Ch.~7]{ec2}. If $f(x_1,\dots,x_k)$ is a
symmetric function, then write $f(M)$ for
$f(\theta_1,\dots,\theta_k)$, where $\theta_1,\dots,\theta_k$ are the
eigenvalues of $M\in U(k)$. The Schur functions $s_\lambda$ for
$\ell(\lambda)\leq k$ are the irreducible characters of $U(k)$, so by
the orthogonality of characters we have for
partitions $\lambda,\mu$ of length at most $k$ that 
   \beq \int_{M\in U(k)} s_\lambda(M)\overline{s_\mu(M)}dM =
      \delta_{\lambda\mu}. \label{eq:chor} \eeq
Now $\mathrm{tr}(M)^n = p_1(M)^n$, where
$p_1(x_1,\dots,x_k)=x_1+\cdots +x_k$. The symmetric function
$p_1(x_1,\dots,x_k)^n$ has 
the expansion \cite[Exam.~1.5.2]{macd}\cite[Cor.~7.12.5]{ec2}
  \[ p_1(x_1,\dots,x_k)^n = \sum_{\substack{\lambda\vdash n\\
      \ell(\lambda)\leq k}} f^\lambda s_\lambda(x_1,\dots,x_k), \]
where $f^\lambda$ is the number of SYT of shape $\lambda$ as in
Section~\ref{sec:enum}. (This formula is best understood algebraically
as a consequence of the Schur-Weyl duality between $\sn$ and $U(k)$
\cite[Ch.~7, Appendix~2]{ec2}, although it can be proved without any
recourse to representation theory.) Hence from
          equation~(\ref{eq:chor}) we obtain 
   \begin{align*} E(|\mathrm{tr}(M)^n|^2) & =  \int_{M\in U(k)}
          p_1(M)^n\overline{p_1(M)}^ndM\\ & = 
     \sum_{\substack{\lambda\vdash n\\ \ell(\lambda)\leq k}} \left(
     f^\lambda\right)^2. \end{align*}
Comparing with equation (\ref{eq:uknsum}) completes the proof.
\end{proof}

Many variations of Theorem~\ref{thm:trmn} have been investigated. For
instance, we can replace $\mathrm{tr}(M)^n$ by more general symmetric
functions of the eigenvalues, such as $\mathrm{tr}(M^m)^n$, or we can
replace $U(k)$ with other classical groups, i.e., $O(k)$ and
$\mathrm{Sp}(2k)$. For further information, see Rains \cite{rains1}. 

\section{Pattern avoidance.} \label{sec:pav}
In this and the following two sections we consider some
generalizations of increasing/de\-creas\-ing subsequences of permutations.
In this section and the next we look at other kinds of subsequences of
permutations, while in Section~\ref{sec:match} we generalize the
permutations themselves.

We have defined $u_k(n)$ to be the number of permutations in $\sn$
with no increasing subsequence of length $k+1$. We can instead
prohibit other types of subsequences of a fixed length, leading to the
currently very active area of \emph{pattern avoidance}.

Given $v=b_1\cdots b_k\in\fs_k$, we say that a permutation
$w=a_1\cdots a_n\in\sn$ \emph{avoids} $v$ if it contains no
subsequence $a_{i_1}\cdots a_{i_k}$ in the same relative order as $v$,
i.e., no subsequence $a_{i_1}\cdots a_{i_k}$ satisfies the condition:
 \[ \forall\ 1\leq r<s\leq k,\ \ a_{i_r}<a_{i_s}
    \Leftrightarrow b_r<b_s. \]
Thus a permutation $w$ satisfies is$(w)<k$ if and only if it is
$12\cdots k$-avoiding, and similarly satisfies ds$(w)<k$ if and only
if it is $k(k-1)\cdots 1$-avoiding. What can be said about the set
$\sn(v)$ of permutations $w\in \sn$ that are $v$-avoiding? In
particular, when are there formulas and recurrences for
$s_n(v):=\#\sn(v)$ similar to those to those of
Theorem~\ref{thm:gessel} and Example~\ref{ex:ukn}? 

The vast subject of pattern avoidance, as a generalization of
avoiding long increasing and decreasing subsequences, began in 1968
with Knuth \cite[Exer.~2.2.1.5]{knuth3}. He showed
in connection with a problem on stack sorting
that $s_n(312)$ is the Catalan number $C_n$. (See also
\cite[Exer.~6.19(ff)]{ec2}.) By obvious symmetries this result,
together with equation (\ref{eq:u2n}), shows that $s_n(v)=C_n$ for all
$v\in\fs_3$. A fundamental paper directly connecting 321-avoiding and
231-avoiding permutations was written by Simion and Schmidt \cite{s-s}.

Wilf first raised the question of investigating $s_n(v)$ for
$v\in\fs_k$ when $k\geq 4$. Here is a brief summary of some highlights
in this burgeoning area. For further information, see e.g.\
\cite[Chs.~4,5]{bona} and the special issue \cite{eljc}. Call two
permutations $u,v\in\fs_k$ \emph{equivalent}, denoted $u\sim v$, if
$s_n(u)=s_n(v)$ for 
all $n$. Then there are exactly three equivalence classes of
permutations $u\in\fs_4$. One class contains 1234, 1243, and 2143 (and
their trivial symmetries), the second contains 3142 and 1342, and the
third 1324. The values of $s_n(1234)$ are given by equation
(\ref{eq:u3n}) and of $s_n(1342)$ are given by the generating function
  \[ \sum_{n\geq 0}s_n(1342)x^n = \frac{32x}{1+20x-8x^2-
    (1-8x)^{3/2}}, \]
a result of B\'ona \cite{bona2}.The enumeration of 1324-avoiding
permutations in $\sn$ remains open.
  
\begin{figure}
\centering
 \centerline{\psfig{figure=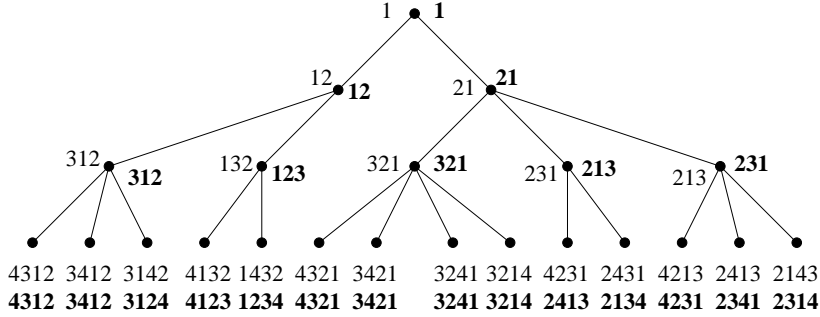}}
\caption{The generating tree for $\sn(123)$ and $\sn(132)$}
\label{fig:gentree}
\end{figure}

Let us mention one useful technique for showing the equivalence of
permutations in $\fs_k$, the method of \emph{generating trees}
introduced by Chung, Graham, Hoggatt and Kleiman \cite{cghk} and
further developed by West \cite{west2}\cite{west}\cite{west3} and
others. Given $u\in\fs_k$, the \emph{generating tree} ${\cal T}_u$ is
the tree with vertex set $\bigcup_{n\geq 1}\sn(u)$, and with $y$ a
descendent of $w$ if $w$ is a subsequence of $y$ (an actual
subsequence, not the pattern of a subsequence). For many pairs
$u,v\in\fs_k$ we have ${\cal T}_u\cong {\cal T}_v$, showing in
particular that $s_n(u)=s_n(v)$ for all $n$, i.e., $u\sim v$. In many
cases in fact the two trees will have no automorphisms, so the
isomorphism ${\cal T}_u\rightarrow {\cal T}_v$ is unique, yielding a
canonical bijection $\sn(u)\rightarrow\sn(v)$. A unique isomorphism
holds for instance when $u=123$ and $v=132$.  Figure~\ref{fig:gentree}
shows the first four levels of the trees ${\cal T}_{123}\cong {\cal
T}_{132}$, labelled by elements of both ${\cal T}_{123}$ and ${\cal
T}_{132}$ (boldface). This tree can also be defined recursively by the
condition that the root has two children, and if vertex $x$ has $k$
children, then the children of $x$ have $2,3,\dots,k+1$ children. For
further information about trees defined in a similar recursive manner,
see Banderier et al.\ \cite{band}.

Given $v\in\fs_k$, let
  \[ F_v(x) = \sum_{n\geq 0}s_n(v)x^n. \]
It is not known whether $F_v(x)$ is always algebraic or even the
weaker condition of being \emph{D-finite}, which is equivalent to
$s_n(v)$ being P-recursive \cite{rs:df}\cite[{\S}6.4]{ec2}. A
long-standing conjecture, known as the \emph{Stanley-Wilf conjecture},
stated that for all $v\in\fs_k$ there is a $c>1$ such that
$s_n(v)<c^n$. In this case
$L_v:=\lim_{n\rightarrow\infty}s_n(v)^{1/n}$ exists and satisfies
$1<L_v<\infty$ \cite{arratia}. For instance, equation (\ref{eq:limuk})
asserts that $L_{12\cdots k}=(k-1)^2$. The Stanley-Wilf conjecture was
proved by Marcus and Tardos \cite{m-t} by a surprisingly simple
argument. It is known that for $k\geq 4$ the permutation $12\cdots k$
neither maximizes nor minimizes $L_v$ for $v\in\fs_k$ \cite{aerw}, but
it is not known which permutations do achieve the maximum or minimum.

An interesting aspect of pattern avoidance was considered by Albert
\cite{albert}. Let $X$ be a finite subset of
$\fs_2\cup\fs_3\cup\cdots$, and let $s_n(X)$ denote the number of
permutations $w\in\sn$ avoiding all permutations $v\in X$. We say that
$X$ is \emph{proper} if it 
doesn't contain both the identity permutation $12\cdots j$ for some
$j$ and the ``reverse identity'' $k\cdots 21$ for some $k$. It is easy
to see (using Theorem~\ref{thm:e-s}) that $X$ is proper if and only if
$s_n(X)>0$ for all $n\geq 1$. For $w\in\sn$ let $\is^X(w)$ be the
length of the longest subsequence of $w$ avoiding all $v\in X$, and
let $E_X(n)$ denote the expectation of $\is^X(w)$ for uniform
$w\in\sn$, so
  $$ E_X(n)=\frac{1}{n!}\sum_{w\in\sn} \is^X(w). $$
Albert then makes the following intriguing conjecture.

\begin{conj} \label{conj:albert}
If $X$ is proper then
  \beq \lim_{n\to\infty} s_n(X)^{1/n} = \frac 14\left(
     \lim_{n\to\infty}\frac{E_X(n)}{\sqrt{n}}\right)^2. 
    \label{eq:albert} \eeq
(The limit on the left-hand side exists by \cite{m-t} and a simple
generalization of \cite{arratia}, while the limit on the right-hand
side is only conjectured to exist.)
\end{conj}

Consider for instance the case $X=\{ 21\}$. Then $s_n(X)=1$ for all
$n\geq 1$, so the left-hand side of \eqref{eq:albert} is 1. On the
other hand, $\is^X(w)= \is(w)$, so the right-hand side is also 1 by
equation~\eqref{eq:vkexp}. More generally, Albert proves
Conjecture~\ref{conj:albert} for $X=\{12\cdots k\}$ and hence (by
symmetry) for $X=\{ k\cdots 21\}$ \cite[Prop.~4]{albert}.

Let us mention that permutations avoiding a certain pattern $v$ or a
finite set of patterns have arisen naturally in a variety of contexts.
For instance, an elementary result of Tenner \cite{tenner} asserts
that the interval $[\mathrm{0},w]$ in the Bruhat order of $\sn$ is a
boolean algebra if and only if $w$ is 321 and 3412-avoiding, and that
the number of such permutations $w$ is the Fibonacci number
$F_{2n-1}$. The Schubert polynomial $\fs_w$ is a single monomial if
and only if $w$ is 132-avoiding \cite[p.~46]{macd2}. All reduced
decompositions of $w\in \sn$ are connected by the Coxeter relations
$s_is_j=s_js_i$ (i.e., $s_is_{i+1}s_i$ does not appear as a factor in
any reduced decomposition of $w$) if and only if $w$ is 321-avoiding
\cite[Thm.~2.1]{b-j-s}.  \emph{Vexillary permutations} may be defined
as those permutations $w$ such that the stable Schubert polynomial
$F_w$ is a single Schur function or equivalently, whose Schubert
polynomial $\fs_w$ is a flag Schur function (or multi-Schur
function). They turn out to be the same as 2143-avoiding permutations
\cite[(1.27)(iii), (7.24)(iii)]{macd2}, first enumerated by West
\cite[Cor.~3.17]{west2}\cite[Cor.~3.11]{west}. Similarly, permutations
$w\in\sn$ for which the Schubert variety $\Omega_w$ in the complete
flag variety GL$(n,\cc)/B$ is smooth are those permutations that are
4231 and 3412-avoiding (implicit in Ryan \cite{ryan}, based on earlier
work of Lakshmibai, Seshadri, and Deodhar, and explicit in Lakshmibai
and Sandhya \cite{la-sa}). The enumeration of such ``smooth
permutations'' in $\sn$ is due to Haiman
\cite{bona3}\cite{haiman}\cite[Exer.~6.47]{ec2}, viz.,
  \[ \sum_{n\geq 0} \sn(4231,3412)x^n = \frac{1}{1-x-\frac{x^2}{1-x} 
        \left( \frac{2x}{1+x-(1-x)C(x)}-1 \right)}, \]
where $C(x)=\sum_{n\geq 0}C_nx^n = (1-\sqrt{1-4x})/2x$.
See Billey and Lakshmibai \cite{b-l} for further information.
As a final more complicated example, Billey and Warrington \cite{b-w}
show that a permutation $w\in\sn$ has a number of nice properties
related to the Kazhdan-Lusztig polynomials $P_{x,w}$ if and only if $w$
avoids 321, 46718235, 46781235, 56718234, and 56781234. These
permutations were later enumerated by Stankova and West \cite{s-w}.
A database of ``natural occurrences'' of pattern avoidance can be
found at a website \cite{tenner:pav} maintained by B. Tenner.

The subjects of pattern avoidance and increasing/decreasing
subsequences can be considered together, by asking for the
distribution of $\is(w)$ or $\ds(w)$ where $w$ ranges over a
pattern-avoiding class $\sn(v)$. (For that matter, one can look at the
distribution of $\is(w)$ or $\ds(w)$ where $w$ ranges over \emph{any}
``interesting'' subset of $\sn$.) For instance, Reifegerste
\cite[Cor.~4.3]{reife} shows that for $k\geq 3$,
  \beq \#\{w\in\sn(231)\st \is(w)<k\} = \frac{1}{n}\sum_{i=1}^{k-1}
          {n\choose i}{n\choose i-1}, \label{eq:reife} \eeq
a sum of \emph{Narayana numbers} \cite[Exer.~6.36]{ec2}. Note that the
left-hand side of (\ref{eq:reife}) can also be written as $\#\{w\in
\sn(231, 12\cdots k)\}$, the number of permutations in $\sn$ avoiding
both 231 and $12\cdots k$. Asymptotic results were obtained by
Deutsch, Hildebrand, and Wilf \cite{d-h-w} for the distribution of
$\is(w)$ when $v= 231, 132$, and 321. Their result for $v=132$ is the
following. 

\begin{theorem} \label{thm:dhw}
For $w\in\sn(132)$ the random variable $\is(w)$ has mean $\sqrt{\pi
  n}+O(n^{1/4})$ and standard deviation
$\sqrt{\pi(\frac{\pi}{3}-1)}\sqrt{n}+O(n^{1/4})$. Moreover, for any
$t>-\sqrt{\pi}$ we have 
   \[ \lim_{n\rightarrow\infty}\mathrm{Prob}\left( \frac{\is(w)
     -\sqrt{\pi n}}{\sqrt{n}}\leq t\right) = \sum_{j\in\zz}
      (1-2j^2(t+\sqrt{\pi})^2)e^{-(t+\sqrt{\pi})^2j^2}. \]
\end{theorem} 

The proof of Theorem~\ref{thm:dhw} is considerably easier than its
counterpart for $w\in\sn$ (Theorem~\ref{thm:bdj}) because there is a
relatively simple formula for the number $f(n,k)$ of permutations
$w\in\sn(132)$ satisfying $\is(w)<k$, viz.,
  \[ f(n,k) = 2\sum_{i=\lceil -n/(k+1)\rceil}^{\lfloor (n+1)/(k+1)
         \rfloor} \left( {2n\choose n+i(k+1)} -\frac 14
      {2n+2\choose n+1+i(k+1)}\right). \]

A number of variations and generalizations of pattern-avoiding
permutations have been investigated. In particular, we can look at
patterns where some of the terms must appear consecutively.  This
concept was introduced by Babson and Steingr\'{\i}msson \cite{b-s} and
further investigated by Claesson \cite{claes} and others.  For
instance, the generalized pattern 1--32 indicates a subsequence $a_i
a_j a_{j+1}$ of a permutation $w=a_1 a_2 \cdots a_n$ such that
$a_i<a_{j+1}<a_j$.  The hyphen in the notation 1--32 means that the
first two terms of the subsequence need not be consecutive. The
permutations in $\fs_4$ avoiding 1--32 are all $C_4=14$ permutations
avoiding 132 (in the previous sense, so avoiding 1--3--2 in the
present context) together with 2413. A typical result, due to Claesson
\cite[Props.~2 and 5]{claes}, asserts that
  \[ \#\sn(\mbox{1--23}) = \#\sn(\mbox{1--32}) = B(n), \]
the number of partitions of the set $[n]$ (a \emph{Bell number}
\cite[{\S}1.4]{ec1}).

\section{Unimodal and alternating subsequences.} \label{sec:alt}
We briefly discuss two variations of increasing/decreasing
subsequences of a different flavor from those considered above. There
is considerable room for further work in this area.

Early work of Chung \cite{chung} and Steele \cite{steele} deals with
$k$-unimodal subsequences. A sequence is $k$-\emph{unimodal} if it is
a concatenation of (at most) $k+1$ monotone sequences. (In more
traditional terminology, a sequence is $k$-unimodal if it has at most
$k+1$ \emph{alternating runs} \cite[{\S}1.2]{bona}.) Thus a
0-unimodal sequence is just an increasing or decreasing sequence, and
every such sequence is $k$-unimodal for all $k$. The sequences 41235
and 24531 are 1-unimodal. Chung showed that every $w\in\sn$ has a
1-unimodal subsequence of length $\left\lceil \sqrt{2n+\frac 14}-\frac
  12\right\rceil$, and that this result is best possible. She
conjectured that if $E_k(n)$ is the expected length of the longest
$k$-unimodal subsequence of a random permutation $w\in \sn$, then
$E_k(n)/\sqrt{n}$ approaches a limit $c_k$ as $n\rightarrow\infty$.
Steele proved this conjecture and showed that $c_k= 2\sqrt{k+1}$ by
deducing it from the monotone ($k=0$) case.

A sequence $b_1 b_2\cdots b_k$ of integers is \emph{alternating} if
  \[ b_1>b_2<b_3>b_4<\cdots b_k. \]
For instance, there are five alternating permutations in $\fs_4$,
viz., 2143, 3142, 4132, 3241, 4231. If $E_n$ denotes the number of
alternating permutations in $\sn$, then a famous result of Andr\'e
\cite{andre}\cite[{\S}3.16]{ec2} states that
  \beq \sum_{n\geq 0} E_n\frac{x^n}{n!} = \sec x + \tan x. 
   \label{eq:engf} \eeq
The numbers $E_n$ were first considered by Euler (using (\ref{eq:engf})
as their definition) and are known as \emph{Euler numbers}. Sometimes
$E_{2n}$ is called a \emph{secant number} and $E_{2n-1}$ a
\emph{tangent number}.   

We can try to extend the main results on increasing/decreasing
subsequences to alternating subsequences. In particular, given
$w\in\sn$ let $\as(w)=\as_n(w)$ denote the length of the longest
alternating subsequence of $w$, and define
  \[ b_k(n) =\#\{w\in\sn\st \as(w)\leq k\}. \]
Thus $b_1(n)=1$ (corresponding to the permutation $12\cdots
n$), $b_k(n)=n!$ if $k\geq n$, and $b_n(n)-b_{n-1}(n)=E_n$. 
Note that we can also define $b_k(n)$ in terms of pattern avoidance,
viz., $b_k(n)$ is the number of $w\in \sn$ avoiding all
$E_{k+1}$ alternating permutations in $\fs_{k+1}$.

Unlike the situation for $u_k(n)$ (defined by (\ref{eq:ukn})), there
are ``nice'' explicit generating functions and formulas for
$b_k(n)$. The basic reason for the existence of such explicit
results is the following (easily proved) key lemma.

\begin{lemma} \label{lemma:alt}
For any $w\in\sn$, there exists an alternating subsequence of $w$ of
maximum length that contains $n$. 
\end{lemma}

From Lemma~\ref{lemma:alt} it is straightforward to derive a
recurrence satisfied by $a_k(n):=b_k(n)-b_{k-1}(n)$, viz.,
  \beq a_k(n) = \sum_{j=1}^n {n-1\choose j-1}
   \sum_{2r+s=k-1}(a_{2r}(j-1)+a_{2r+1}(j-1))a_s(n-j). 
   \label{eq:aknrec} \eeq
Now define
  \[ B(x,t) = \sum_{k,n\geq 0}b_k(n)t^k\frac{x^n}{n!}. \]
It follows  from the recurrence (\ref{eq:aknrec}) (after
some work \cite{rs:alt}) that
  \beq B(x,t)=\frac{1+\rho+2te^{\rho x}+(1-\rho)
     e^{2\rho x}}{1+\rho-t^2+(1-\rho-t^2)e^{2\rho x}}, \label{eq:fxt} 
  \eeq
where $\rho=\sqrt{1-t^2}$. Alternatively (as pointed out by M.
B\'ona), let $G(n,k)$ denote the number of $w\in \sn$ with $k$
alternating runs as defined at the beginning of this section. Then
equation~(\ref{eq:fxt}) is a consequence of the relation $a_k(n)=
\frac 12(G(n,k-1)+G(n,k))$ and known facts about $G(n,k)$ summarized
in \cite[{\S}1.2]{bona}.

It can be deduced from equation~(\ref{eq:fxt}) (shown with assistance
from I. Gessel) that
  \[ b_k(n) = \frac{1}{2^{k-1}} \sum_{\substack{i+2j\leq k\\ i\equiv
     k\,(\mathrm{mod}\,2)}} 
     (-2)^j{k-j\choose (k+i)/2}{n\choose j}i^n. \]
For instance,
  \[ b_1(n)  =  1,\ \  b_2(n) =  2^{n-1},\ \
        b_3(n) =  \frac 14(3^n-2n+3), \ \
        b_4(n) =  \frac 18(4^n-(2n-4)2^n).
  \]
From equation (\ref{eq:fxt}) it is also easy to compute the moments  
  \[ M_k(n) = \frac{1}{n!}\sum_{w\in\sn} \as_n(w)^k. \]
For instance,
  \begin{align*} \sum_{n\geq 1}M_1(n)x^n & =  
       \left.\frac{\partial B(x,t)}{\partial t}\right|_{t=1}\\
    & =  \frac{6x-3x^2+x^3}{6(1-x)^2}, \end{align*}
from which we obtain 
  \beq M_1(n) = \left\{ \begin{array}{rl} 1, & n=1\\[.05in]
      \displaystyle\frac{4n+1}{6}, & n>1. \end{array} \right. 
   \label{eq:m1n} \eeq
Similarly the variance of $\as_n(w)$ is given by 
  \beq \mathrm{var}(\as_n) = \frac{8}{45}n-\frac{13}{180},\ n\geq 4. 
    \label{eq:varasn} \eeq
It is surprising that there are such simple explicit formulas, in
contrast to the situation for $\is(w)$ (equation (\ref{eq:en})). 

It is natural to ask for the limiting distribution of $\as_n$,
analogous to Theorem~\ref{thm:bdj} for $\is_n$. The following result
was shown independently by R. Pemantle \cite{pemantle} and H. Widom
\cite{widom}. It can also be obtained by showing that the polynomials
$\sum_k a_k(n)t^k$ have (interlacing) real zeros, a consequence of the 
connection between $a_k(n)$ and $G(n,k)$ mentioned above and a result
of Wilf \cite[Thm.~1.41]{bona}.

\begin{theorem}
We have for random (uniform) $w\in\sn$ and all $t\in\rr$ that
  $$ \lim_{n\to\infty} \mathrm{Prob}\left(
   \frac{\as_n(w)-2n/3}{\sqrt{n}} \leq t\right) = G(t), $$
where $G(t)$ is Gaussian with variance $8/45$:
  $$ G(t) = \frac{1}{\sqrt{\pi}} \int_{-\infty}^{t\sqrt{45}/4}
        e^{-s^2}ds. $$
\end{theorem}
\indent For further information on longest alternating subsequences,
see the paper \cite{rs:alt}.

\section{Matchings.} \label{sec:match}
The subject of pattern containment and avoidance discussed in
Section~\ref{sec:pav} provides one means to extend the
concept of increasing/decreasing subsequences of permutations. In this
section we will consider a different approach, in which permutations
are replaced with other combinatorial objects.  We will be concerned
mainly with (complete) \emph{matchings} on $[2n]$, which may be
defined as partitions $M=\{B_1,\dots,B_n\}$ of $[2n]$ into $n$
two-element blocks $B_i$. Thus $B_1\cup B_2\cup\cdots \cup B_n=[2n]$,
$B_i\cap B_j=\emptyset$ if $i\neq j$, and $\#B_i=2$.  (These
conditions are not all independent.) Alternatively, we can regard a
matching $M$ as a fixed-point free involution $w_M$ of $[2n]$, viz.,
if $B_i=\{a,b\}$ then $w_M(a)=b$. We already considered increasing and
decreasing subsequences of fixed-point free involutions in
Section~\ref{sec:sym}. In that situation, however, there is no
symmetry interchanging increasing subsequences with decreasing
subsequences. Here we consider two alternative statistics on 
matchings (one of which is equivalent to decreasing subsequences)
which have the desired symmetry.

\begin{figure}
\centering
 \centerline{\psfig{figure=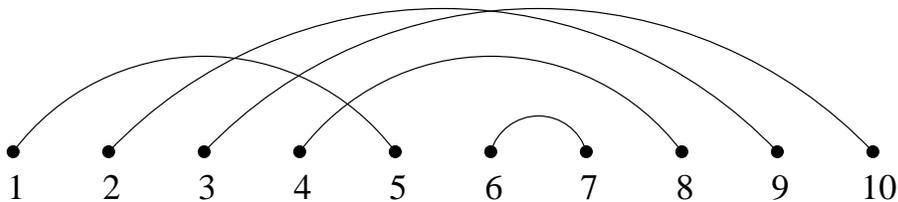}}
\caption{A matching on $[10]$}
\label{fig:matchex}
\end{figure}

Write $\fm_n$ for the set of matchings on $[2n]$. We represent a
matching $M\in\fm_n$ by a diagram of $2n$ vertices $1,2,\dots,2n$ on a
horizontal line in the plane, with an arc between vertices $i$ and $j$
and lying above the vertices if $\{i,j\}$ is a block of $M$.
Figure~\ref{fig:matchex} shows the diagram corresponding to the
matching
  \[ M = \{\{1,5\},\{2,9\},\{3,10\},\{4,8\},\{6,7\}\}. \]

Let $M\in\fm_n$. A \emph{crossing} of $M$ consists of two
arcs $\{i,j\}$ and $\{k,l\}$ with $i<k<j<l$. Similarly a
\emph{nesting} of $M$ consists of two arcs $\{i,j\}$ and $\{k,l\}$
with $i<k<l<j$. The maximum number of mutually crossing arcs of $M$ is
called the $\emph{crossing number}$ of $M$, denoted $\cro(M)$.
Similarly the \emph{nesting number} $\nes(M)$ is the maximum number of
mutually nesting arcs. For the matching $M$ of
Figure~\ref{fig:matchex}, we have $\cro(M)=3$ (corresponding to the
arcs $\{1,5\}$, $\{2,9\}$), and $\{3,10\}$), while also $\nes(M)=3$
(corresponding to $\{2,9\}$, $\{4,8\}$, and $\{6,7\}$).  

It is easy to see that $\ds(w_M)=2\cdot\nes(M)$, where $w_M$ is the
fixed-point free involution corresponding to $M$ as defined above.
However, it is not so clear whether $\cro(M)$ is connected with
increasing/decreasing subsequences. To this end, define 
  \[ f_n(i,j) =\#\{ M\in\fm_n\st \cro(M)=i,\ \nes(M)=j\}. \]
It is well-known that
  \beq \sum_j f_n(0,j) = \sum_i f_n(i,0) = C_n. \label{eq:nocr} \eeq
In other words, the number of matchings $M\in\fm_n$ with no crossings
(or with no nestings) is the Catalan number $C_n$. For crossings this
result goes back to Errera \cite{err}\cite[Exer.~6.19(n,o)]{ec2}; for
nestings see \cite{rs:cat}. Equation (\ref{eq:nocr}) was given the
following generalization by Chen et al. \cite{cddsy}.

\begin{theorem} \label{thm:fnij}
For all $i,j,n$ we have $f_n(i,j)=f_n(j,i)$.
\end{theorem}

\begin{figure}
\centering
\centerline{\psfig{figure=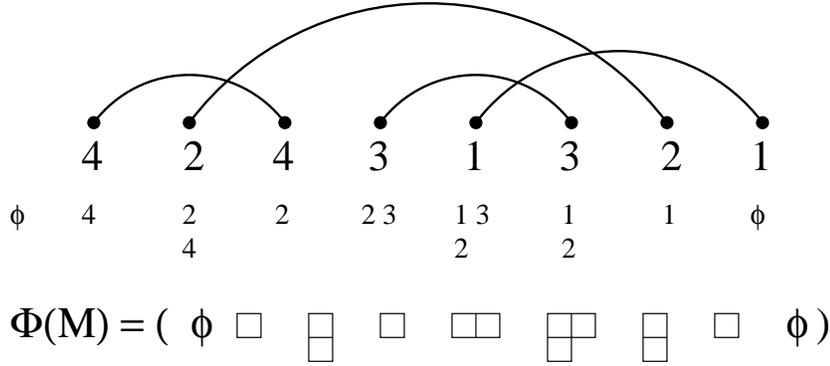}}
\caption{A correspondence between matchings and oscillating tableaux}
\label{fig:oscrsk}
\end{figure}

Theorem~\ref{thm:fnij} is proved by using a version of RSK first
defined by the author (unpublished) and then extended by Sundaram
\cite{sundaram}. Define an \emph{oscillating tableau} of shape
$\lambda\vdash n$ and length $k$ to be a sequence
  \[ \emptyset =\lambda^0, \lambda^1,\dots,\lambda^k=\lambda \]
of partitions $\lambda^i$ such that (the diagram of) $\lambda^{i+1}$
is obtained from $\lambda^i$ by either adding or removing a square.
(Note that if we add a square each time, so $k=n$, then we obtain an
SYT of shape $\lambda$.) Oscillating tableaux were first defined
(though not with that name) by Berele \cite{berele} in connection with
the representation theory of the symplectic group. Given a matching
$M\in\fm_n$, define an oscillating tableau
$\Phi(M)=(\lambda^0,\lambda^1,\dots,\lambda^{2n})$ of length $2n$ and
shape $\emptyset$ as follows. Label the right-hand endpoints of the
arcs of $M$ by $1,2,\dots,n$ from right-to-left. Label each left-hand
endpoint with the same label as the right-hand endpoint. Begin with
the empty tableau $T_0=\emptyset$. Let $a_1,\dots, a_{2n}$ be the
sequence of labels, from left-to-right. Once $T_{i-1}$ has been
obtained, define $T_i$ to be the tableau obtained by row-inserting $a_i$
into $T_{i-1}$ (as defined in Section~\ref{sec:enum}) if $a_i$ is the
label of a left-hand endpoint of an 
arc; otherwise $T_i$ is the tableau obtained by deleting $a_i$ from
$T_{i-1}$. Let $\lambda^i$ be the shape of $T_i$, and set
  \[ \Phi(M) = (\emptyset=\lambda^0,
  \lambda^1,\dots,\lambda^{2n}=\emptyset). \]
See Figure~\ref{fig:oscrsk} for an example. It is easy to see that
$\Phi(M)$ is an oscillating tableau of length $2n$ 
and shape $\emptyset$. With a little more work it can be shown that in
fact the map $M\mapsto \Phi(M)$ is a \emph{bijection} from $\fm_n$ to the
set ${\cal O}_n$ of all oscillating tableaux of length $2n$ and shape
$\emptyset$. As a consequence we have the enumerative formula
  \beq \#{\cal O}_n = (2n-1)!! := 1\cdot 3\cdot 5\cdots(2n-1), 
     \label{eq:enma} \eeq
the number of matchings on $[2n]$. The key fact about the
correspondence $\Phi$ for proving Theorem~\ref{thm:fnij} is the
following ``oscillating analogue'' of Schensted's theorem
(Theorem~\ref{thm:sch}). 

\begin{theorem} \label{thm:osch}
Let $M\in\fm_n$, and suppose that $\Phi(M)=(\lambda^0,\lambda^1,\dots,
\lambda^{2n})$. Then $\nes(M)$ is equal to the most number of columns
of any $\lambda^i$, while $\cro(M)$ is equal to the most number of
rows of any $\lambda^i$.
\end{theorem}

Theorem~\ref{thm:fnij} is an easy consequence of
Theorem~\ref{thm:osch}. For let $\Phi(M)'$ be the oscillating tableau
obtained by conjugating all the partitions in $\Phi(M)$, and let $M'$
be the matching satisfying $\Phi(M')=\Phi(M)'$. By
Theorem~\ref{thm:osch} we have $\cro(M)=\nes(M')$ and
$\nes(M)=\cro(M')$. Since the map $M\mapsto \Phi(M)$ is a bijection
we have that the operation $M\mapsto M'$ is a bijection from
$\fm_n$ to itself that interchanges cr and ne, and the proof follows. 

The above argument and equation~\eqref{eq:wr} show that the operation
$M\mapsto M'$ on matchings is a natural analogue of the operation of
reversal on permutations. Unlike the case of permutations, we don't
know a simple ``direct'' operation on matchings that interchanges cr
with ne.

Let $\tilde{f}^\lambda_n$ denote the number
of oscillating tableaux of shape $\lambda$ and length $n$, so
$\tilde{f}^\lambda_n=0$ unless $n\equiv |\lambda|\,
(\mathrm{mod}\,2)$. A generalization for $\tilde{f}^\lambda_n$ of the
hook-length formula (equation (\ref{eq:frt})) is due to Sundaram
\cite[Lemma~2.2]{sundaram}, viz.,
  \[ \tilde{f}^\lambda_n = {n\choose k} (n-k-1)!!\,f^\lambda, \]
where $\lambda\vdash k$ and where of course
$f^\lambda$ is evaluated by the usual hook-length formula
(\ref{eq:frt}). (Set $(-1)!!=1$ when $n=k$.) Note that an oscillating
tableau $(\lambda^0, 
\lambda^1,\dots, \lambda^{2n})$ of shape $\emptyset$ may be regarded
as a pair $(P,Q)$ of oscillating tableaux of the same shape
$\lambda=\lambda^n$ and length $n$, viz.,
  \begin{align*} P & =  (\lambda^0,\lambda^1,\dots,\lambda^n)\\
        Q & =  (\lambda^{2n},\lambda^{2n-1},\dots,\lambda^n). \end{align*}
Hence we obtain the following restatement of equation (\ref{eq:enma}):
  \beq \sum_\lambda \left( \tilde{f}_n^\lambda\right)^2 = (2n-1)!!, 
     \label{eq:tfn2} \eeq
where $\lambda$ ranges over all partitions. (The partitions $\lambda$
indexing a nonzero summand are those satisfying $\lambda\vdash k\leq
n$ and $k\equiv n\,(\mathrm{mod}\,2)$.)

Equation (\ref{eq:tfn2}) suggests, in analogy to equations
(\ref{eq:fla}) and (\ref{eq:sumfl2}), a connection between
$\tilde{f}_n^\lambda$ and representation theory. Indeed, there is a
$\cc$-algebra $\mathfrak{B}_n(x)$, where $x$ is a real parameter,
which is semisimple for all but finitely many $x$ (and such that these
exceptional $x$ are all integers) and which has a basis that is
indexed in a natural way by matchings $M\in\fm_n$. In particular,
$\dim \mathfrak{B}_n(x)=(2n-1)!!$. This algebra was first defined
by Brauer \cite{brauer} and shown to be the centralizer algebra of the
action of the orthogonal group O$(V)$ on $V^{\otimes n}$ (the $n$th
tensor power of $V$), where $\dim V=k$ and $x=k$. It is also the
centralizer algebra of the action of the symplectic group Sp$(2k)$ on
$V^{\otimes n}$, where now $\dim V=2k$ and $x=-2k$. When
$\mathfrak{B}_n(x)$ is semisimple, its irreducible representations
have dimension $\tilde{f}_n^\lambda$, so we obtain a
representation-theoretic explanation of equation~(\ref{eq:tfn2}). For
further information see e.g.\ Barcelo and Ram \cite[App.~B6]{ba-ra} 

Because $\ds(w_M)=2\cdot\nes(M)$ and because $\cro_n$ and $\nes_n$
have the same distribution by Theorem~\ref{thm:fnij}, the asymptotic
distribution of $\cro_n$ and $\nes_n$ on $\fm_n$ reduces to that of
$\ds_{2n}$ on $\mfi_{2n}^*$, which is given by
Theorem~\ref{thm:sym}(b). We therefore obtain the following result.

\begin{theorem} \label{thm:crdis}
We have for random (uniform) $M\in\fm_n$ and all $t\in\rr$ that 
  \[ \lim_{n\rightarrow \infty} \mathrm{Prob}\left(
    \frac{\nes_n(M)-\sqrt{2n}}{(2n)^{1/6}}\leq \frac t2\right) = 
    F(t)^{1/2}\exp\left( \frac 12\int_t^\infty u(s)ds\right). \]
The same result holds with $\nes_n$ replaced with $\cro_n$. 
\end{theorem}

We can also consider the effect of bounding \emph{both}
$\cro(M)$ and $\nes(M)$ as $n\rightarrow \infty$. The analogous
problem for $\is(w)$ and $\ds(w)$ is not interesting, since
by Theorem~\ref{thm:e-s} there are no permutations $w\in\sn$
satisfying $\is(w)\leq p$ and $\ds(w)\leq q$ as soon as
$n>pq$. Let 
  \begin{align*} h_{p,q}(n) & =  \#\{M\in\fm_n\st \cro(M)\leq p,\
         \nes(M)\leq q\}\\
     H_{p,q}(x) & =  \sum_{n\geq 0}h_{p,q}(n)x^n. \end{align*}
It follows from the bijection $\Phi:\fm_n\rightarrow {\cal
O}_n$ and a simple application of the transfer-matrix method
\cite[{\S}4.7]{ec1} that $H_{p,q}(x)$ is a rational function of $x$
\cite[{\S}5]{cddsy}. For instance,  
  \[ H_{1,1}(x) =  \frac{1}{1-x},\ \
     H_{1,2}(x) =  \frac{1-x}{1-2x},\ \ 
     H_{1,3}(x) =  \frac{1-2x}{1-3x+x^2} \]
  \[ H_{2,2}(x) =  \frac{1-5x+2x^2}{(1-x)(1-5x)},\ \
     H_{2,3}(x) =  \frac{1-11x+30x^2-23x^3+4x^4}
            {(1-x)(1-3x)(1-8x+4x^2)} \]
  \[ H_{3,3}(x) =  \frac{1-24x+186x^2-567x^3+690x^4-285x^5+15x^6}
    {(1-x)(1-19x+83x^2-x^3) (1-5x+6x^2-x^3)^2}. \]
Christian Krattenthaler pointed out that $h_{p,q}(n)$ can be
interpreted as counting certain walks in an alcove of the affine Weyl
group $\tilde{C}_n$. It then follows from a result of Grabiner
\cite[(23)]{grab} that all reciprocal zeros of the denominator of
$H_{p,q}(x)$ are of the form
  \[ 2(\cos(\pi r_1/m)+\cdots+\cos(\pi r_j/m)), \]
where each $r_i\in\zz$ and $m=p+q+1$. All these numbers for fixed $m$
belong to an extension of $\qq$ of degree $\phi(2m)/2$, where $\phi$
is the Euler phi-function. As a consequence, every irreducible factor
(over $\qq$) of the denominator of $H_{p,q}(x)$ has degree dividing
$\phi(2m)/2$.

Theorem~\ref{thm:fnij} can be extended to objects other than
matchings, in particular, arbitrary set partitions. (Recall that we
have defined a matching to be a partition of $[2n]$ into $n$ 2-element
blocks.) In this situation oscillating tableaux are replaced by
certain sequences of partitions known as \emph{vacillating
tableaux}. See \cite{cddsy} for further details. Vacillating tableaux
were introduced implicitly (e.g., \cite[(2.23)]{h-r}) in connection
with the representation theory of the \emph{partition algebra}
$\mathfrak{P}_n$, a semisimple algebra whose dimension is the Bell
number $B(n)$. See Halverson and Ram \cite{h-r} for a survey of the
partition algebra. Vacillating tableaux and their combinatorial
properties were made more explicit by Chen, et al.\ \cite{cddsy} and
by Halverson and Lewandowski \cite{h-l}. An alternative approach based
on ``growth diagrams'' to vacillating tableaux and their nesting and
matchings was given by Krattenthaler \cite{kratt}.  It remains open to
find an analogue of Theorem~\ref{thm:crdis} for the distribution of
$\cro(\pi)$ or $\nes(\pi)$ (as defined in \cite{cddsy}) for arbitrary
set partitions $\pi$.


\end{document}